\documentstyle[12pt,prl,aps,graphicx]{revtex}

\newcommand{\be}{\begin{equation}}
\newcommand{\ee}{\end{equation}}
\newcommand{\np}{\newpage}
\begin{document} \draft \date{\today} 
\title{Permutations and Primes}
\author{Zakir F. Seidov\\
\it {Department of Physics, Ben-Gurion University of the Negev,
Beer-Sheva 84105, Israel}} \maketitle \begin{abstract}
The problem of digital sets (DS) {\it all} permutations of which
generate primes
is discussed. Such sets may include only four types of digits: 1, 3, 7 and
9. The direct
computations show that for N-digit (ND) integers such DS's 
occur at N = 1, 2, 3, and are absent in the 4 - 10 interval of N. On the other
hand in N = 19, 23, 317 and 1031 cases (as
well as in N = 2 case) 
the formal "total" permutation is provided by  "repunits", integers with
all digits being 1. 
  The existence/nonexistence of other (not
{repunits}) full-permutation DS's for arbitrary large N is an open  question with
probable negative answer. Maximal-permutation DS's with maximal number of
primes are given for various numbers of digits. 
\end{abstract}

The primes are real pearls among all integers and are subject
of  long-standing (several thousand year) researches.
The many amazing results are established on primes [1, 2, 3].\\ 
The number of primes
are infinite (by the way the largest known prime with 2 098 960
digits, $2^{6\, 972\, 593}-1$, was found quite recently [1])
 and their general behavior  with $n$ is  at present 
well known and  is described, e.g., by some built-in functions
of MATHEMATICA [4]: PrimeQ[n] gives True if n is prime and False
otherwise,
Prime[n] gives nth prime (with assumption Prime[1] = 2), and 
PrimePi[n] gives
number of all primes $\leq\,n$.\\
Using mainly MATHEMATICA's means I discuss here a 
problem of the permutations of given digital set all giving prime.  
Consider some
examples first. Already for two-digit primes one may note pairs
13/31, 17/71, 37/73, 79/97 - all primes,
(one may add  here fifth "pair" 11/11 as well),  while 19 "loses" its
partner, 91 which is not a prime, see Table 2.
One may wonder: 
 is it possible, in the  more-digit cases,
that full permutations of some DS all give 
primes?\\
In 3D case (Table 2, case B) the answer is positive, there are 3 full
permutation DS's: 113, 199 and 337, that is,
for instance, all three possible permutations 113, 131 and 331 give
primes.
Other two interesting DS's are 179 and 379 each giving even more primes,
 4, which is however less than a corresponding total number of possible
permutations - 6. Another DS 
 with 6 permutations, 137, gives only 3 primes.\\

Notice that in this paper the "basic" DS's are written as first
   pemutations in lexicografic order and they themselves
are not    necessarily primes. Two smallest such sets are
   119 and 133 which are not primes while two other 
   members of each permutation family (191, 911 and
   313, 331, respectively) are primes. (And the largest such DS
considered in this note is $7_{2}9_8\,=\,7799999999$, see Table 8.)\\
   Here is the place to clarify the point of terminology: speaking about, e.g.
137 as DS, I consider it as what in MATHEMATICA is $List\{1, 3, 7\}$
with 3 elements;
at another hand, 137 is also used as 3-digit integer  in our usual 10-base
arithmetic system. To avoid cumbersome notations I assume this
not-too-rigorous approach. \\
It is evident that {\it full} permutation may give {\it all}
primes 
for DS's comprising {\it only} of digits 1, 3, 7 and 9, and {\it none} of
0, 2, 4, 5,
6, 8 digits.
This greatly reduces the 
field of search, removing all integers containing at least one of
the digits 0, 2, 4, 5, 6 or 8 .\\ 
The program was written in MATHEMATICA for searching such 
full permutations among "1-3-7-9" primes with the negative
answer for 4 to 10-digit integers. In all cases number of primes provided
by any given DS is {\it less} than number of {\it all} possible
permutations, see attached Figures and Tables.
This is the main (and negative ) result, but  some other remarks may be 
made on the by-products of calculations.
  
  1. Primes among all "1-3-7-9" integers my be briefly described by
"N-P-D" code as follows: \\
 1 - 2 - 2, 2 - 10 - 6,
3 - 30 - 12, 4 - 63 - 14, 5 - 249 - 35, 6 - 757 - 54, 7 - 2 709 - 74,
8 - 9 177 - 101,
9 - 33 191 - 142, 10 - 118 912 - 184,\\  where N is number of digits, P is
number of primes and D is number of DS's.\\
  For  all primes (with arbitrary digits)  "N-P-D" figures are:\\ 
1 - 4 - 4, 2 - 21 - 17, 3 - 143 - 86, 4 - 1 061 - 336,
 5 - 8 363 - 1 109, 6 - 68 906 - 2 967, 7 - 586 081 - 7 041, 
8 - 5 096 876 - 15 259, 9 - 45 086 079 - ?, 10 - 404 204 977 - ?. \\
There is of some interest that for smaller N "mean productivity" (P/D) 
of "1-3-7-9" DS's are larger than that of all DS's, while starting 
from N=5 number of primes per DS is larger for all DS:
\begin{verbatim}
 Table 1. Number of N digit primes per basic DS as function of N 
 for "1-3-7-9" DS's and all DS's
 ========================
  N            P/D
 ========================
        "1379"      all
 ========================
  1     1.0000     1.0000
  2     1.6667     1.2353
  3     2.5000     1.6628
  4     4.5000     3.1577
  5     7.1143     7.5410
  6    14.0185    23.2241
  7    36.6081    83.2383
  8    90.8614   334.0242
  9   233.7394      ?
 10   646.2609      ?
\end{verbatim}

I do not present      
values of D and distribution of number of primes per given DS in general case by two
reasons, first due to necessity of the
too time-consuming calculations, and second 
due to ambiguity of the problem, e.g. what to do with primes with zeros,
as in this case some permutation of corresponding DS, with zeros
at the beginning, give integers (and possible primes) with less
number digits. First three such primes are 101, 103, 107 and last
three such 10 digit primes  are $9_{7}701, \ \ 9_{7}703\ \   \mbox {and} \
\ \ 9_{7}707$.

   Another observation is that while there is a very
   strong correlation between p and c, the more rich 
   permutation family does not necessarily 
   give more primes; cf. 139 vs 113, or 1139 vs 1379.

3. As among primes with 4 to 10 digits there is no single full permutaion, 
it is interesting to look for {\it maximal}  permutations. 
Quite surprisingly, these "maximal" DS's
are not at all among  "1-3-7-9" DS's. From Table 2
one may note, that maximal "1-3-7-9" DS's with 4 digits are
1379 with 7 primes and 1139 with 8 primes. At the same time, 
there are 4 "ordinary" 
DS's (1349, 1457, 3479, 3679) with 9 primes, two (1579 and 1789)
with 10, and two "super-sets" 1237 and 1279 with 11 primes. 
To be able to give more primes DS should comprise the even digits
as well, not only odd ones.
This fact is really very suprising, because number distribution of digits
in, for example, 4D primes is not quite at all random, with great favor to
 digits 1, 3, 7, 9:\\
Count[n] = 217, 603, 359, 602, 326, 327, 336, 574, 321, 579, at n=0 to
9,\\
each digit from "1-3-7-9" family is roughly two times more common than
any one of other family "0-2-4-5-6-8". Inspite of this, DS's of "mixed
races" are more productive in making primes. 
A very interesting observation which may find its application in the
fields very far from primes and integers.\\
4. Some maximal permutation "1379" DS's for the larger numbers of digits
are:
\begin{verbatim}
                 b     c     p    c/p
               ======================= 
      5D       11339   15    30  .5000
               13379   18    60  .3000
               13799   29    60  .4833

      6D      113779   60   180  .3333
              133379   35   120  .2917
              133799   55   180  .3056

      7D     1113799  113   420  .2690
             1133779  182   630  .2889
             1137799  169   630  .2683
    
      8D    11333779  419  1680  .2494
            11337779  403  1680  .2399
            11377999  397  1680  .2363

      9D   113337799 1388  7560  .1836
           111337799 1525  7560  .2017
           113377999 1550  7560  .2050

     10D  1113337799 4555 25200  .1808
          1133777999 4606 25200  .1828
          1133377799 4384 25200  .1740
\end{verbatim}
In all cases c[i]$<<$p[i], and relation c/p has a tendency of decreasing
with increasing N, that is in some sense, "probability" of appearing
of full permutation among "1379" primes diminishes with increasing 
number of digits. One may considered it as a hint
that the full-permutation "1379" DS is absent for arbitrary
large N.\\     
  I conclude with guess that the existence/nonexistence of (not
{repunits}) DS's for arbitrary large N is an open  question with
probable negative answer.

\begin{figure} \includegraphics[scale=.6]{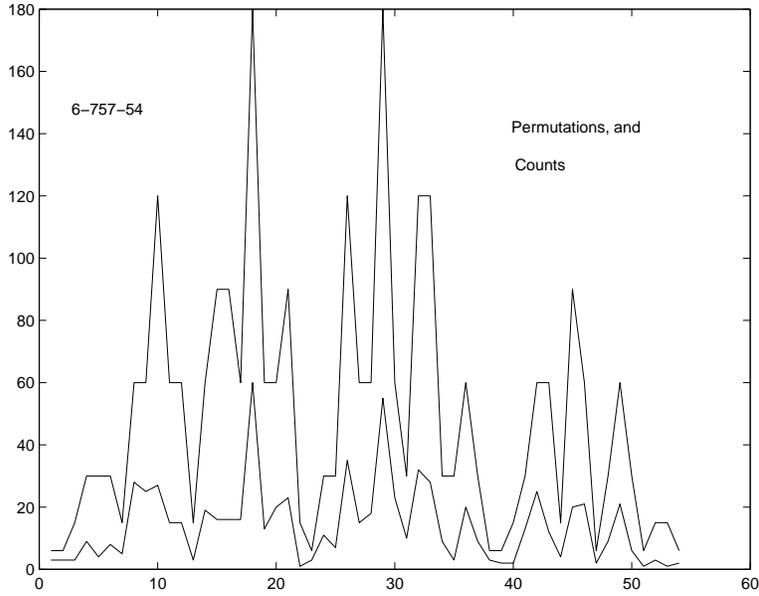} 
\caption{6 digit case. Abscissae - order number of basic digit
sets (BDS); see Table 4. Ordinates -
upper curve: number of permutations of BDS as function of order number of
BDS;
lower curve: number of primes given by BDS as function of order number of
BDS.}
\end{figure}
\begin{figure} \includegraphics[scale=.8]{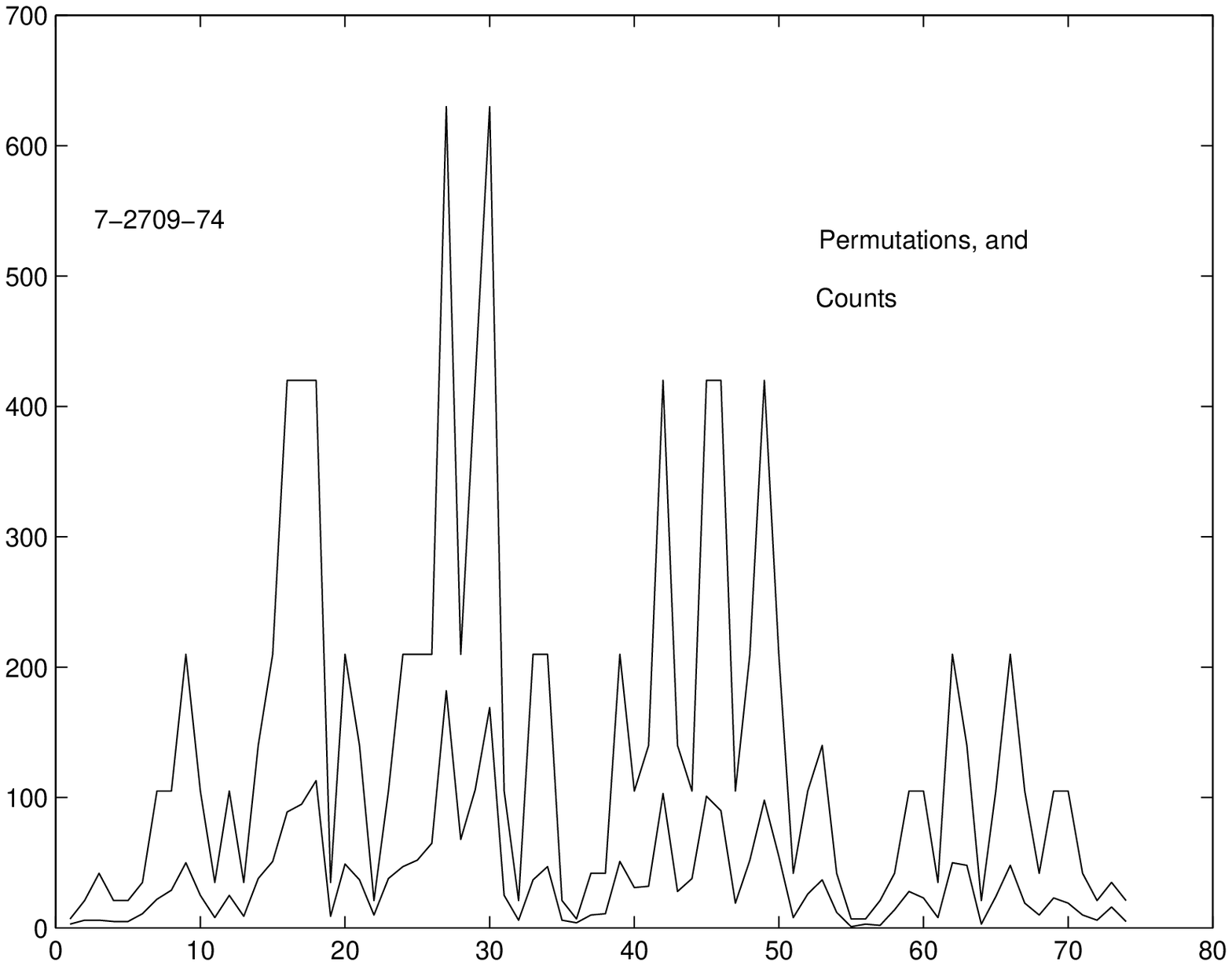}
\caption{7 digit case. Abscissae - order number of basic digit
sets (BDS); see Table 5. Ordinates -
upper curve: number of permutations of BDS as function of order number of
BDS;
lower curve: number of primes given by BDS as function of order number of
BDS.} \end{figure}
 \begin{figure} \includegraphics[scale=.8]{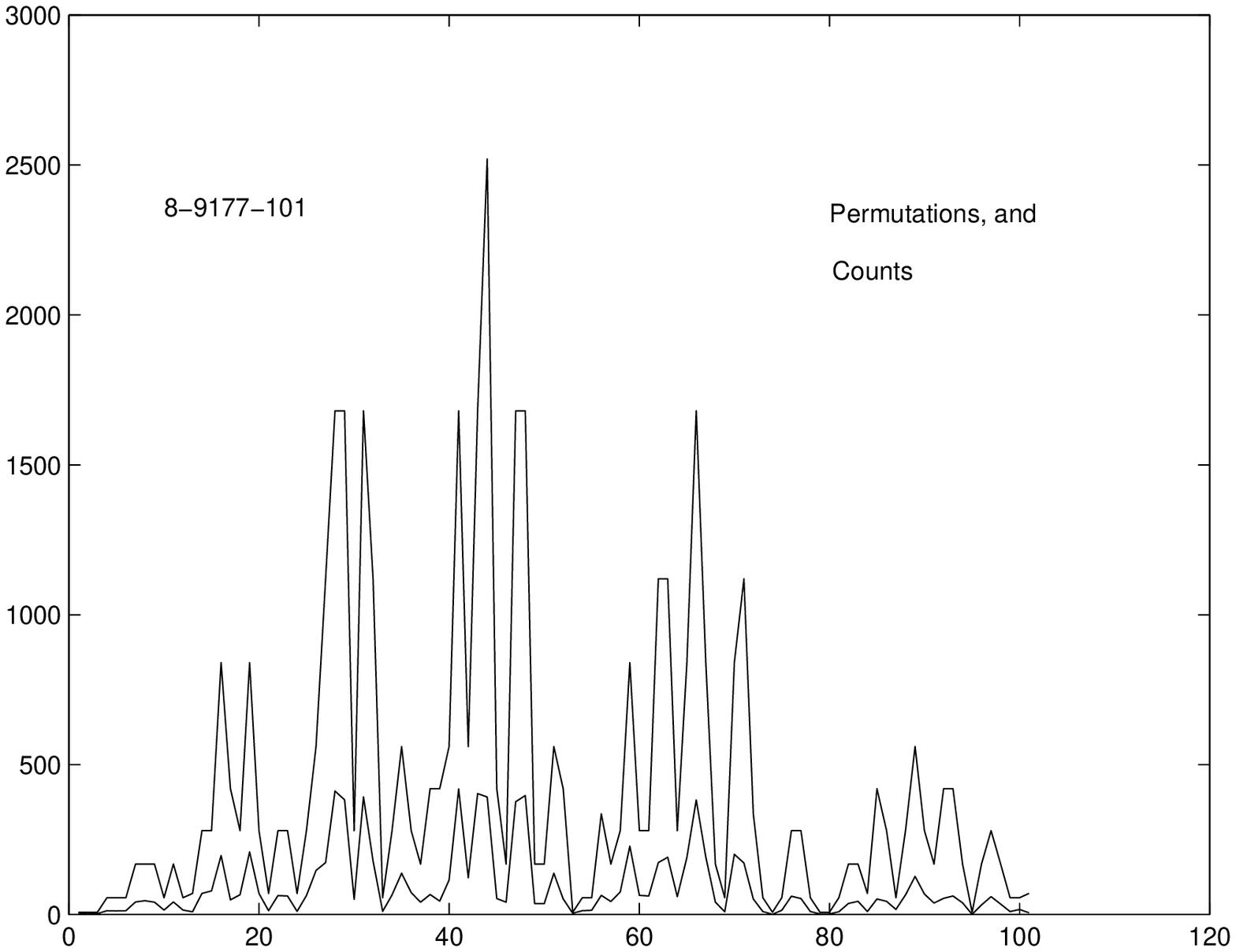} 
\caption{8 digit case. Abscissae - order number of basic digit
sets (BDS); see Table 6. Ordinates -
upper curve: number of permutations of BDS as function of order number of
BDS;
lower curve: number of primes given by BDS as function of order number of
BDS.}\end{figure}
\begin{figure} \includegraphics[scale=.8]{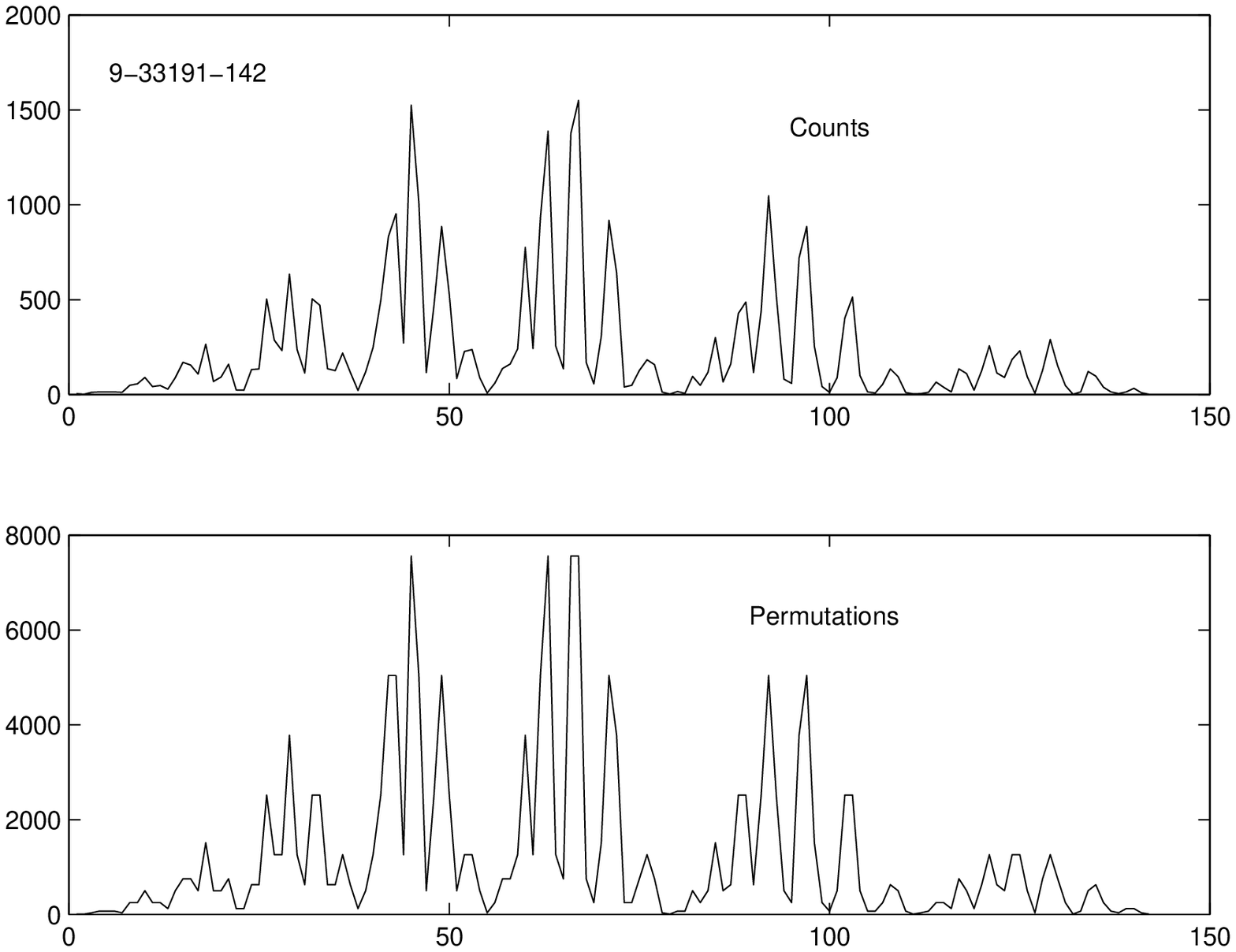} 
\caption{9 digit case. Abscissae - order number of basic digit
sets (BDS); see Table 7. Ordinates -
upper panel: number of primes given by BDS as function of order number of
lower panel: number of permutations of BDS as function of order number of
BDS}\end{figure}
\begin{figure} \includegraphics[scale=.8] {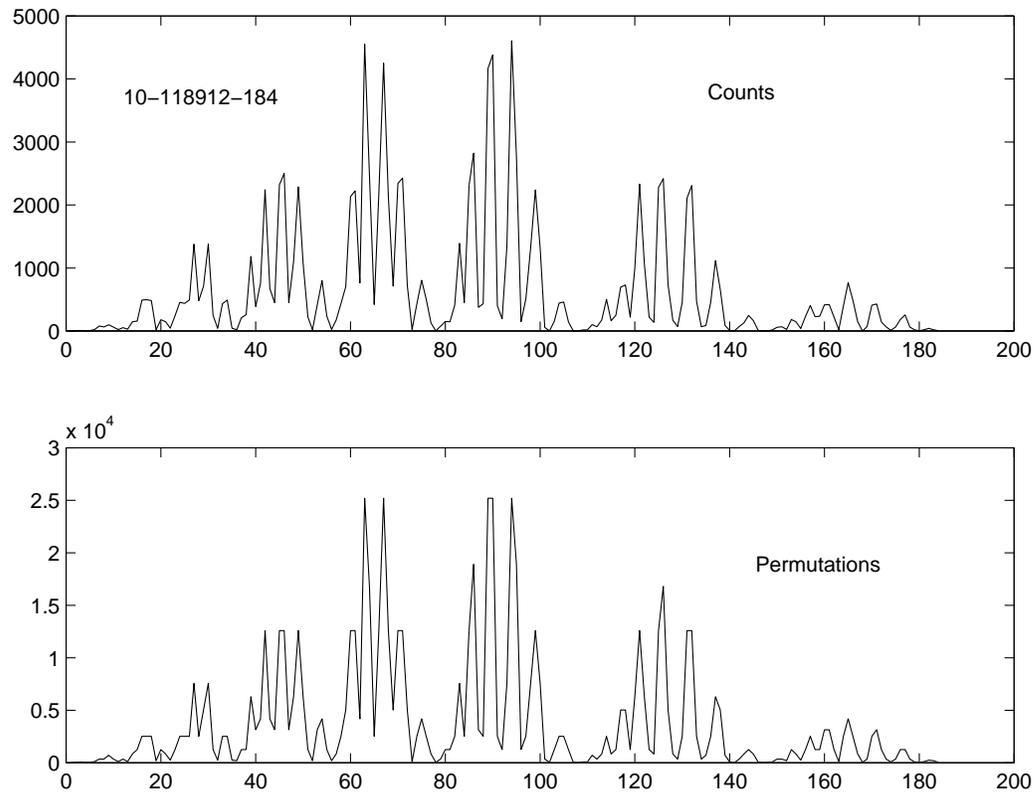} 
\caption{ 10 digit case. Abscissae - order number of basic digit
sets (BDS); see Table 8. Ordinates -
upper panel: number of primes given by BDS as function of order number of
lower panel: number of permutations of BDS as function of order number of
BDS}\end{figure}
\np
\begin{verbatim}
  Table 2. A: 2-10-6 case. b - basic digit sets, 
   c - counts, p - permutations.    
   Total number 10 = Sum[c[i]].
          B: 3-30-12 case.
   c[i] = p[i] at 3 cases, for b=113, 199 and 337. 
          C: 4-63-14 case.   
   c[i] < p[i] at all cases.
   Note that basic digit sets are written as first
   ones in lexicografic order and they are not 
   necessary primes themselves. Two least such sets are
   119 and 133 which are not primes while two other 
   members of each permutation family (191, 911 and
   313, 331, respectively) are primes.
   Another observation is that while there is a very
   strong correlation between p and c, the more rich 
   permutation family does not necessarily 
   give more primes; cf. 139 vs 113, or 1139 vs 1379
    
   =============================================== 
       A                B                C  
   ================================================    
    b    c  p        b     c   p       b     c   p
   ================================================
   11   1  1        113   3   3       1117   2   4        
   13   2  2        119   2   3       1139   8  12
   17   2  2        133   2   3       1333   2   4 
   19   1  2        137   3   6       1337   5  12
   37   2  2        139   2   6       1339   5  12 
   79   2  2        179   4   6       1379   7  24
   ----------       199   3   3       1399   6  12
       10           337   3   3       1777   3   4
                    377   1   3       1799   5  12
                    379   4   6       1999   2   4
                    779   2   3       3337   3   4
                    799   1   3       3379   6  12 
                    -----------       3779   5  12  
                         30           3799   4  12
                                      -----------
                                            63
\end{verbatim} \np \begin{verbatim}


  Table 3. 5-249-35 case. b - 35 basic digit sets, 
  c - counts, p - permutations.    
  Total number 249 = Sum[c[i]].
   c[i] < p[i]  at all i. 

    b    c  p        b     c   p      b     c   p
   ==============================================
  11113   3   5    11117   2   5    11119   1   5       
  11137   5  20    11177   5  10    11179   9  20
  11333   4  10    11339  15  30    11377  10  30 
  11399  13  30    11777   3  10    11779  12  30 
  11999   5  10    13333   2   5    13337   9  20   
  13339   9  20    13379  18  60    13399   8  30 
  13777   8  20    13799  29  60    13999   7  20 
  17777   1   5    17779   5  20    17999   7  20
  19999   1   5    33377   3  10    33379   7  20 
  33779  10  30    33799   9  30    37777   2   5 
  37799  10  30    37999   9  20    77779   4   5 
  77999   3  10    79999   1   5   
  -----------------------------------------------
                                          249
\end{verbatim}
\np \begin{verbatim}
                  Table 4. 6-757-54 Case
                    ====================
                       b     c        p
                    ====================
                    111113   3        6
                    111119   3        6
                    111133   3        15
                    111137   9        30
                    111139   4        30
                    111179   8        30
                    111199   5        15
                    111337   28       60
                    111377   25       60
                    111379   27       120
                    111779   15       60
                    111799   15       60
                    113333   3        15
                    113339   19       60
                    113377   16       90
                    113399   16       90
                    113777   16       60
                    113779   60       180
                    113999   13       60
                    117779   20       60
                    117799   23       90
                    119999   1        15
                    133333   3        6
                    133337   11       30
                    133339   7        30
                    133379   35       120
                    133399   15       60
                    133777   18       60
                    133799   55       180
                    133999   23       60
                    137777   10       30
                    137779   32       120
                    137999   28       120
                    139999   9        30
                    177779   3        30
                    177799   20       60
                    179999   9        30
                    199999   3        6
                    333337   2        6
                    333377   2        15
                    333379   13       30
                    333779   25       60
                    333799   12       60
                    337777   4        15
                    337799   20       90
                    337999   21       60
                    377777   2        6
                    377779   9        30
                    377999   21       60
                    379999   6        30
                    777779   1        6
                    777799   3        15
                    779999   1        15
                    799999   2        6
 \end{verbatim}

\np \begin{verbatim}
                    Table 5. 7-2709-74 Case 
                    =======================
                      b       c         p
                    =======================
                    1111117   3         7
                    1111133   6         21
                    1111139   6         42
                    1111177   5         21
                    1111199   5         21
                    1111333   11        35
                    1111337   22        105
                    1111339   29        105
                    1111379   50        210
                    1111399   25        105
                    1111777   8         35
                    1111799   25        105
                    1111999   9         35
                    1113337   38        140
                    1113377   51        210
                    1113379   89        420
                    1113779   95        420
                    1113799   113       420
                    1117777   9         35
                    1117799   49        210
                    1117999   37        140
                    1133333   10        21
                    1133339   38        105
                    1133377   47        210
                    1133399   52        210
                    1133777   65        210
                    1133779   182       630
                    1133999   68        210
                    1137779   106       420
                    1137799   169       630
                    1139999   25        105
                    1177777   6         21
                    1177799   37        210
                    1177999   47        210
                    1199999   6         21
                    1333333   4         7
                    1333337   10        42
                    1333339   11        42
                    1333379   51        210
                    1333399   31        105
                    1333777   32        140
                    1333799   103       420
                    1333999   28        140
                    1337777   38        105
                    1337779   101       420
                    1337999   90        420
                    1339999   19        105
                    1377779   52        210
                    1377799   98        420
                    1379999   55        210
                    1399999   8         42
                    1777799   26        105
                    1777999   37        140
                    1799999   12        42
                    1999999   1         7
                    3333337   3         7
                    3333377   2         21
                    3333379   14        42
                    3333779   28        105
                    3333799   23        105
                    3337777   8         35
                    3337799   50        210
                    3337999   48        140
                    3377777   3         21
                    3377779   24        105
                    3377999   48        210
                    3379999   19        105
                    3777779   10        42
                    3777799   23        105
                    3779999   19        105
                    3799999   10        42
                    7777799   6         21
                    7777999   16        35
                    7799999   5         21
\end{verbatim}

\np \begin{verbatim}
                    Table 6. 8-9177-101 Case
                    ==========================
                      b         c         p
                    ==========================             
                    11111113   4          8
                    11111117   3          8
                    11111119   3          8
                    11111137   13         56
                    11111179   12         56
                    11111333   13         56
                    11111339   42         168
                    11111377   46         168
                    11111399   41         168
                    11111777   15         56
                    11111779   42         168
                    11111999   15         56
                    11113333   9          70
                    11113337   70         280
                    11113339   79         280
                    11113379   197        840
                    11113399   49         420
                    11113777   66         280
                    11113799   209        840
                    11113999   70         280
                    11117777   13         70
                    11117779   63         280
                    11117999   62         280
                    11119999   11         70
                    11133337   63         280
                    11133377   147        560
                    11133379   173        1120
                    11133779   412        1680
                    11133799   383        1680
                    11137777   51         280
                    11137799   392        1680
                    11137999   178        1120
                    11177777   11         56
                    11177779   65         280
                    11177999   138        560
                    11179999   73         280
                    11333339   41         168
                    11333377   67         420
                    11333399   45         420
                    11333777   115        560
                    11333779   419        1680
                    11333999   123        560
                    11337779   403        1680
                    11337799   392        2520
                    11339999   53         420
                    11377777   41         168
                    11377799   376        1680
                    11377999   397        1680
                    11399999   36         168
                    11777779   36         168
                    11777999   138        560
                    11779999   52         420
                    13333333   4          8
                    13333337   13         56
                    13333339   14         56
                    13333379   64         336
                    13333399   43         168
                    13333777   75         280
                    13333799   228        840
                    13333999   64         280
                    13337777   62         280
                    13337779   173        1120
                    13337999   191        1120
                    13339999   60         280
                    13377779   190        840
                    13377799   382        1680
                    13379999   192        840
                    13399999   41         168
                    13777777   9          56
                    13777799   201        840
                    13777999   172        1120
                    13799999   52         336
                    13999999   10         56
                    17777777   1          8
                    17777779   14         56
                    17777999   61         280
                    17779999   53         280
                    17999999   10         56
                    19999999   2          8
                    33333337   2          8
                    33333379   9          56
                    33333779   37         168
                    33333799   44         168
                    33337777   10         70
                    33337799   52         420
                    33337999   44         280
                    33377777   17         56
                    33377779   65         280
                    33377999   127        560
                    33379999   68         280
                    33777779   38         168
                    33777799   54         420
                    33779999   62         420
                    33799999   39         168
                    37777777   1          8
                    37777799   32         168
                    37777999   60         280
                    37799999   35         168
                    37999999   10         56
                    77777999   17         56
                    77779999   6          70
 \end{verbatim}

\np \begin{verbatim}
                    
                       Table 7. 9-33191-142 Case
                    ===============================
                        b         c         p
                    ===============================
                    111111113   4           9
                    111111119   1           9
                    111111133   11          36
                    111111137   14          72
                    111111139   14          72
                    111111179   14          72
                    111111199   11          36
                    111111337   48          252
                    111111377   56          252
                    111111379   90          504
                    111111779   41          252
                    111111799   47          252
                    111113333   28          126
                    111113339   89          504
                    111113377   169         756
                    111113399   155         756
                    111113777   108         504
                    111113779   265         1512
                    111113999   69          504
                    111117779   92          504
                    111117799   159         756
                    111119999   24          126
                    111133333   22          126
                    111133337   131         630
                    111133339   135         630
                    111133379   503         2520
                    111133399   286         1260
                    111133777   231         1260
                    111133799   634         3780
                    111133999   237         1260
                    111137777   113         630
                    111137779   504         2520
                    111137999   470         2520
                    111139999   136         630
                    111177779   126         630
                    111177799   218         1260
                    111179999   117         630
                    111199999   21          126
                    111333337   119         504
                    111333377   248         1260
                    111333379   495         2520
                    111333779   833         5040
                    111333799   954         5040
                    111337777   272         1260
                    111337799   1525        7560
                    111337999   1012        5040
                    111377777   116         504
                    111377779   470         2520
                    111377999   886         5040
                    111379999   528         2520
                    111777779   84          504
                    111777799   227         1260
                    111779999   237         1260
                    111799999   89          504
                    113333333   8           36
                    113333339   59          252
                    113333377   137         756
                    113333399   161         756
                    113333777   240         1260
                    113333779   776         3780
                    113333999   242         1260
                    113337779   931         5040
                    113337799   1388        7560
                    113339999   255         1260
                    113377777   136         756
                    113377799   1376        7560
                    113377999   1550        7560
                    113399999   170         756
                    113777777   56          252
                    113777779   308         1512
                    113777999   918         5040
                    113779999   642         3780
                    113999999   39          252
                    117777779   48          252
                    117777799   127         756
                    117779999   183         1260
                    117799999   157         756
                    119999999   12          36
                    133333333   2           9
                    133333337   16          72
                    133333339   6           72
                    133333379   95          504
                    133333399   48          252
                    133333777   117         504
                    133333799   300         1512
                    133333999   66          504
                    133337777   161         630
                    133337779   428         2520
                    133337999   487         2520
                    133339999   116         630
                    133377779   440         2520
                    133377799   1047        5040
                    133379999   522         2520
                    133399999   81          504
                    133777777   58          252
                    133777799   719         3780
                    133777999   885         5040
                    133799999   254         1512
                    133999999   42          252
                    137777777   11          72
                    137777779   89          504
                    137777999   403         2520
                    137779999   513         2520
                    137999999   100         504
                    139999999   13          72
                    177777779   8           72
                    177777799   54          252
                    177779999   135         630
                    177799999   94          504
                    179999999   10          72
                    333333337   3           9
                    333333377   5           36
                    333333379   11          72
                    333333779   65          252
                    333333799   38          252
                    333337777   14          126
                    333337799   135         756
                    333337999   110         504
                    333377777   23          126
                    333377779   126         630
                    333377999   257         1260
                    333379999   113         630
                    333777779   90          504
                    333777799   185         1260
                    333779999   230         1260
                    333799999   93          504
                    337777777   7           36
                    337777799   128         756
                    337777999   290         1260
                    337799999   150         756
                    337999999   47          252
                    377777777   1           9
                    377777779   14          72
                    377777999   121         504
                    377779999   96          630
                    377999999   38          252
                    379999999   14          72
                    777777799   5           36
                    777779999   13          126
                    777799999   32          126
                    779999999   9           36
                    799999999   1           9

\end{verbatim}

\np 
\begin{verbatim}
                      Table 8. 10-118912-184 Case
                     ==============================
                         b           c            p
                     ==============================
                     1111111117      2            10
                     1111111133      5            45
                     1111111139      9            90
                     1111111177      4            45
                     1111111199      5            45
                     1111111333     22           120
                     1111111337     79           360
                     1111111339     64           360
                     1111111379     97           720
                     1111111399     62           360
                     1111111777     23           120
                     1111111799     54           360
                     1111111999     29           120
                     1111113337    147           840
                     1111113377    157          1260
                     1111113379    492          2520
                     1111113779    496          2520
                     1111113799    482          2520
                     1111117777     20           210
                     1111117799    180          1260
                     1111117999    147           840
                     1111133333     45           252
                     1111133339    241          1260
                     1111133377    455          2520
                     1111133399    437          2520
                     1111133777    491          2520
                     1111133779   1379          7560
                     1111133999    479          2520
                     1111137779    709          5040
                     1111137799   1379          7560
                     1111139999    246          1260
                     1111177777     41           252
                     1111177799    433          2520
                     1111177999    490          2520
                     1111199999     45           252
                     1111333333     20           210
                     1111333337    210          1260
                     1111333339    259          1260
                     1111333379   1184          6300
                     1111333399    388          3150
                     1111333777    754          4200
                     1111333799   2242         12600
                     1111333999    673          4200
                     1111337777    449          3150
                     1111337779   2320         12600
                     1111337999   2500         12600
                     1111339999    451          3150
                     1111377779   1096          6300
                     1111377799   2288         12600
                     1111379999   1077          6300
                     1111399999    222          1260
                     1111777777     20           210
                     1111777799    414          3150
                     1111777999    799          4200
                     1111799999    235          1260
                     1111999999     24           210
                     1113333337    174           840
                     1113333377    435          2520
                     1113333379    702          5040
                     1113333779   2132         12600
                     1113333799   2225         12600
                     1113337777    763          4200
                     1113337799   4555         25200
                     1113337999   2473         16800
                     1113377777    421          2520
                     1113377779   2171         12600
                     1113377999   4252         25200
                     1113379999   2131         12600
                     1113777779    713          5040
                     1113777799   2343         12600
                     1113779999   2426         12600
                     1113799999    721          5040
                     1117777777     24           120
                     1117777799    428          2520
                     1117777999    806          4200
                     1117799999    503          2520
                     1117999999    127           840
                     1133333333      2            45
                     1133333339     76           360
                     1133333377    151          1260
                     1133333399    145          1260
                     1133333777    410          2520
                     1133333779   1392          7560
                     1133333999    453          2520
                     1133337779   2311         12600
                     1133337799   2823         18900
                     1133339999    376          3150
                     1133377777    431          2520
                     1133377799   4166         25200
                     1133377999   4384         25200
                     1133399999    400          2520
                     1133777777    194          1260
                     1133777779   1346          7560
                     1133777999   4606         25200
                     1133779999   2755         18900
                     1133999999    148          1260
                     1137777779    505          2520
                     1137777799   1317          7560
                     1137779999   2240         12600
                     1137799999   1317          7560
                     1139999999     56           360
                     1177777777      2            45
                     1177777799    153          1260
                     1177777999    442          2520
                     1177799999    459          2520
                     1177999999    142          1260
                     1199999999      4            45
                     1333333333      1            10
                     1333333337     18            90
                     1333333339     16            90
                     1333333379    101           720
                     1333333399     67           360
                     1333333777    174           840
                     1333333799    502          2520
                     1333333999    163           840
                     1333337777    250          1260
                     1333337779    697          5040
                     1333337999    729          5040
                     1333339999    219          1260
                     1333377779    993          6300
                     1333377799   2332         12600
                     1333379999   1066          6300
                     1333399999    218          1260
                     1333777777    137           840
                     1333777799   2279         12600
                     1333777999   2419         16800
                     1333799999    733          5040
                     1333999999    172           840
                     1337777777     69           360
                     1337777779    455          2520
                     1337777999   2107         12600
                     1337779999   2310         12600
                     1337999999    473          2520
                     1339999999     67           360
                     1377777779     87           720
                     1377777799    454          2520
                     1377779999   1118          6300
                     1377799999    659          5040
                     1379999999     86           720
                     1399999999     12            90
                     1777777777      2            10
                     1777777799     70           360
                     1777777999    130           840
                     1777799999    247          1260
                     1777999999    166           840
                     1799999999     11            90
                     1999999999      1            10
                     3333333377      2            45
                     3333333379     16            90
                     3333333779     58           360
                     3333333799     68           360
                     3333337777     25           210
                     3333337799    181          1260
                     3333337999    145           840
                     3333377777     38           252
                     3333377779    233          1260
                     3333377999    405          2520
                     3333379999    226          1260
                     3333777779    233          1260
                     3333777799    415          3150
                     3333779999    419          3150
                     3333799999    221          1260
                     3337777777     21           120
                     3337777799    395          2520
                     3337777999    765          4200
                     3337799999    486          2520
                     3337999999    147           840
                     3377777777      1            45
                     3377777779     70           360
                     3377777999    409          2520
                     3377779999    430          3150
                     3377999999    144          1260
                     3379999999     65           360
                     3777777779     12            90
                     3777777799     62           360
                     3777779999    182          1260
                     3777799999    256          1260
                     3779999999     62           360
                     3799999999     20            90
                     7777777799      2            45
                     7777777999     19           120
                     7777799999     41           252
                     7777999999     21           210
                     7799999999      5            45
 \end{verbatim}


\end{document}